\documentclass{mathincs}
\usepackage{graphicx}
\usepackage{amsfonts}
\usepackage{amsmath,amssymb,bbm,url}

\newtheorem{theorem}{Theorem}
\newtheorem{proposition}[theorem]{Proposition}
\theoremstyle{definition}
\newtheorem{definition}[theorem]{Definition}
\theoremstyle{remark}
\newtheorem{remark}[theorem]{Remark}

\newcommand{\bv}{\mathbf{b}}
\newcommand{\C}{\mathcal{C}}
\newcommand{\Cv}{\mathbf{C}}

\newcommand{\E}{\mathcal{E}}
\newcommand{\fv}{\mathbf{f}}
\newcommand{\Fv}{\mathbf{F}}
\newcommand{\gv}{\mathbf{g}}
\newcommand{\hv}{\mathbf{h}}
\newcommand{\je}[1]{J_{#1}\pi}

\newcommand{\R}{\mathcal{R}}
\newcommand{\RR}{\mathbbm{R}}
\newcommand{\Rc}[1]{\R_{#1}}

\newcommand{\sv}{\mathbf{s}}
\newcommand{\T}{\mathcal{T}}
\newcommand{\uv}{\mathbf{u}}

\newcommand{\U}{\mathcal{U}}
\newcommand{\V}{\mathcal{V}}
\newcommand{\xv}{\mathbf{x}}
\newcommand{\yv}{\mathbf{y}}

\DeclareMathOperator{\im}{im}

\begin{document}
\title{On the Numerical Analysis and Visualisation of Implicit Ordinary
  Differential Equations}
\author{Elishan Braun}
\address{Departimento di Matematica e Fisica, Universit\`a degli Studi
  Roma Tre, Largo San Leonardo Murialdo 1,\\ 00146 Rome, Italy}
\email{elishan@hotmail.de}
\author{Werner M. Seiler}
\address{Institut f\"{u}r Mathematik, Universit\"at Kassel,
Heinrich-Plett-Stra\ss e 40, 34132 Kassel, Germany}
\email{seiler@mathematik.uni-kassel.de}
\author{Matthias Sei\ss}
\address{Institut f\"{u}r Mathematik, Universit\"at Kassel,
Heinrich-Plett-Stra\ss e 40, 34132 Kassel, Germany}
\email{mseiss@mathematik.uni-kassel.de}

\begin{abstract}
  We discuss how the geometric theory of differential equations can be used
  for the numerical integration and visualisation of implicit ordinary
  differential equations, in particular around singularities of the
  equation.  The Vessiot theory automatically transforms an implicit
  differential equation into a vector field distribution on a manifold and
  thus reduces its analysis to standard problems in dynamical systems
  theory like the integration of a vector field and the determination of
  invariant manifolds.  For the visualisation of low-dimensional situations
  we adapt the streamlines algorithm of Jobard and Lefer to 2.5 and 3
  dimensions.  A concrete implementation in \textsc{Matlab} is discussed
  and some concrete examples are presented.
\end{abstract}
\thanks{This work was partially performed as part of the European
  H2020-FETOPEN-2016-2017-CSA project $SC^{2}$ (712689) and partially
  supported by the bilateral project ANR-17-CE40-0036 and DFG-391322026
  SYMBIONT.} 
\maketitle

\section{Introduction}

Most textbooks on the theoretical or the numerical analysis of ordinary
differential equations assume that the equations are given in the solved
form $\uv^{(q)}=\fv(t,\uv,\dot{\uv},\dots,\uv^{(q-1)})$.  Fully implicit
differential equations $\Fv(t,\uv,\dot{\uv},\dots,\uv^{(q)})=0$ exhibit a
much wider range of behaviours including the appearance of
singularities. Even basic questions of the existence and uniqueness of
solutions are much more involved for them and near singularities basically
all numerical methods break down or become at least ill-conditioned.  In
this article, we demonstrate how a geometric approach allows to translate
many of these questions into standard problems for dynamical systems on
manifolds.

Our approach to implicit equations is based on the Vessiot distribution
associated to any differential equation (see \cite{wms:vessconn2,wms:invol}
and references therein) and was already employed in
\cite{wms:aims,wms:singbif}.  In the case of a not underdetermined ordinary
differential equation (which we will exclusively consider in this article),
the Vessiot distribution is almost everywhere one-dimensional and can thus
be locally represented by a vector field.  One-dimensional integral curves
of this vector field correspond to generalised solutions which around
regular points are nothing but the prolongation of classical solutions.
Irregular singularities become stationary points and the local solution
behaviour around them is determined by their invariant manifolds.

In this article, we will describe the ideas underlying a suite of
\textsc{Matlab} routines (developed as part of the first author's master
thesis \cite{eb:nvide}) that allow (i) for the automated numerical
determination and integration of the Vessiot distribution wherever it is
one-dimensional, (ii) for the visualisation of the streamlines of the
Vessiot distribution in 2, 2.5 and 3 dimensions using a method proposed by
Jobard and Lefer \cite{jl:stream} and (iii) for the numerical determination
of invariant manifolds and the reduced dynamics on them via an approach
developed by Beyn and Kle\ss~\cite{bk:invman} and later improved by Eirola
and von Pfaler \cite{ep:teim}.

The article is structured as follows.  The next section recalls the basic
ingredients of the geometric theory of ordinary differential equations: jet
bundles, Vessiot spaces, singularities and generalised solutions.  Section
\ref{sec:numint} is concerned with the numerical integration away from
irregular singularities which boils down to integrating a vector field on a
manifold.  The analysis of the local solution behaviour around irregular
singularities using invariant manifolds of stationary points is the topic
of Section \ref{sec:irreg}.  The following section describes our ansatz for
visualising streamlines in various dimensions.  In Section
\ref{sec:examples}, two concrete examples of a fully nonlinear first-order
and of a quasi-linear second-order equation are treated.  The reader may
find it useful to refer from time to time to this section, as any concept
or construction discussed in earlier sections will be explicitly
demonstrated there.  Finally, some conclusions are given.

\section{Geometry of Differential Equations}
\label{sec:geode}

Differential equations are geometrically modelled via jet bundles
\cite{vvl:homo,pom:eins,sau:jet,wms:invol}.  Let $\pi:\E\rightarrow\T$ be a
fibred manifold with $\dim{\T}=1$ for ordinary differential equations,
e.\,g.\ $\T=\RR$ and $\E=\T\times\RR^{m}$ with $\pi$ the canonical
projection on $\T$.  For simplicity, we work in local coordinates, although
we use throughout a ``global notation''.  As coordinate on the base space
$\T$ we use $t$ and fibre coordinates in the total space $\E$ will be
$\uv=(u^1,\ldots,u^m)$.  The first derivative of $u^{\alpha}$ will be
denoted by $\dot{u}^{\alpha}$; higher derivatives are written in the form
$u^\alpha_k = d^{k}u^\alpha/dt^{k}$.  Adding all derivatives $u^\alpha_k$
with $k\leq q$ (collectively denoted by $\uv_{(q)}$) defines a coordinate
system for the $q$-th order jet bundle $\je{q}$.  There are natural
fibrations $\pi^{q}_{r}:\je{q}\rightarrow\je{r}$ for $r<q$ and
$\pi^{q}:\je{q}\rightarrow\T$ ``forgetting'' all higher derivatives.
Sections $\sigma:\T\rightarrow\E$ of the fibration $\pi$ correspond to
functions $\uv=\sv(t)$, as locally they can always be written in the form
of a graph $\sigma(t)=\bigl(t,\sv(t)\bigr)$.  To such a section $\sigma$,
we associate its \emph{prolongation} $j_q\sigma:\T\rightarrow\je{q}$, a
section of the fibration $\pi^q$ given by
$j_q\sigma(t)=\bigl(t,\sv(t), \dot{\sv}(t),\ddot{\sv}(t),\dots\bigr)$.

The geometry of the $q$-th order jet bundle $\je{q}$ is to a large extent
determined by its \emph{contact structure} describing intrinsically the
relationship between the different types of coordinates.  The \emph{contact
  distribution} is the smallest distribution $\C_{q}\subset T(\je{q})$ that
contains the tangent spaces $T(\im{j_q\sigma})$ of all $q$ times prolonged
sections and any field in it is a \emph{contact vector field}.  In local
coordinates, $\C_{q}$ is generated by one transversal and $m$ vertical
fields (with respect to $\pi^{q}$):
\begin{subequations}\label{eq:ciq}
  \begin{align}
      C_{\mathrm{trans}}^{(q)} &= 
          \partial_t + \sum_{j=0}^{q-1}u^\alpha_{j+1}\partial_{u^\alpha_j}\;,
          \label{eq:ciq1}\\
      C^{(q)}_\alpha &= \partial_{u^\alpha_q}\;,\qquad 
          1\leq\alpha\leq m\;.\label{eq:ciq2}
  \end{align}
\end{subequations}

\begin{proposition}\label{prop:contact}
  A section $\gamma:\T\rightarrow\je{q}$ is of the form
  $\gamma=j_{q}\sigma$ with $\sigma:\T\rightarrow\E$, if and only if
  $T_{\gamma(t)}(\im{\gamma})\subseteq\C_{q}|_{\gamma(t)}$ for all points
  $t\in\T$ where $\gamma$ is defined.
\end{proposition}

Compared with the usual intrinsic geometric definition of a differential
equation, ours allows for certain types of singularities, as it imposes
considerably weaker conditions on the restricted projection $\hat\pi^{q}$
which in the standard definition must be a surjective submersion.  Note
that we do not distinguish between scalar equations and systems.

\begin{definition}\label{def:diffeq}
  An \emph{(ordinary) differential equation} of order $q$ is a submanifold
  $\Rc{q}\subseteq\je{q}$ such that the restriction $\hat\pi^{q}$ of the
  projection $\pi^{q}:\je{q}\rightarrow\T$ to $\Rc{q}$ has a dense image in
  $\T$.  A \emph{(strong) solution} is a (local) section
  $\sigma:\T\rightarrow\E$ such that $\im{j_{q}\sigma}\subseteq\Rc{q}$.
\end{definition}

Locally, a differential equation $\Rc{q}\subseteq\je{q}$ can be described
as the zero set of some smooth functions $\Phi:\je{q}\rightarrow\RR$ which
brings us back to the usual picture of a differential equation.  Checking
whether a function $\sv(t)$ is a solution by entering it and its
derivatives into $\Phi$ corresponds to verifying that
$\im{j_{q}\sigma}\subseteq\Rc{q}$ for the section $\sigma$ defined by
$\sv(t)$.  In the sequel, we will always assume that we are dealing with a
formally integrable equation, i.\,e.\ that there are no hidden
integrability conditions.  We refer to \cite{wms:invol} for more details on
the meaning and the effective verification of this assumption.
Furthermore, we will assume throughout that all considered differential
equations are not underdetermined, i.\,e.\ their general solution depends
only on a finite number of constants.

A key insight of Cartan was to study \emph{infinitesimal solutions} or
\emph{integral elements} of a differential equation
$\Rc{q}\subseteq\je{q}$, i.\,e.\ to consider at any point $\rho\in\Rc{q}$
those linear subspaces $\U_\rho\subseteq T_\rho\Rc{q}$ which are
potentially the tangent space at $\rho$ of a prolonged solution through
$\rho$.  We will follow here an approach pioneered by Vessiot
\cite{ves:int} which is based on vector fields and dual to the more popular
Cartan-K\"ahler theory of exterior differential systems (see
\cite{df:diss,wms:vessconn2,wms:invol} for modern presentations).  By
Proposition \ref{prop:contact}, the tangent spaces $T_\rho(\im{j_q\sigma})$
of prolonged sections at points $\rho\in\je{q}$ are always subspaces of the
contact distribution $\C_q|_{\rho}$.  If the section $\sigma$ is a solution
of $\Rc{q}$, it furthermore satisfies $\im{j_q\sigma}\subseteq\Rc{q}$ by
Definition \ref{def:diffeq} and hence
$T_{\rho}(\im{j_q\sigma})\subseteq T_{\rho}\Rc{q}$ for any point
$\rho\in\im{j_q\sigma}$.  These considerations motivate the following
construction.

\begin{definition}\label{def:vessiot}
  The \emph{Vessiot space} of the differential equation
  $\Rc{q}\subseteq\je{q}$ at the point $\rho\in\Rc{q}$ is that part of the
  contact distribution that lies tangential to $\Rc{q}$, i.\,e.\ the vector
  space
  \begin{equation}\label{eq:vessiot}
    \V_{\rho}[\Rc{q}]=T_{\rho}\Rc{q}\cap\C_q|_{\rho}\;.
  \end{equation}
  The family of all Vessiot spaces of a differential equation $\Rc{q}$ is
  its \emph{Vessiot distribution} $\V[\Rc{q}]$.
\end{definition}

Our discussion of the explicit construction of the Vessiot spaces in the
next section will show that they are almost everywhere one-dimensional
(because of our restriction to not underdetermined equations) and that they
define on an open subset of $\Rc{q}$ a smooth regular distribution.
Following the terminology of Arnold \cite{via:geoode}, we will distinguish
the points on $\Rc{q}$ according to the properties of their Vessiot spaces.

\begin{definition}\label{def:sing}
  A point $\rho\in\Rc{q}$ is \emph{regular}, if $\dim{\V_{\rho}[\Rc{q}]}=1$
  and $\V_{\rho}[\Rc{q}]$ is transversal  relative to the fibration
  $\pi^{q}$.  If $\dim{\V_{\rho}[\Rc{q}]}=1$, but $\V_{\rho}[\Rc{q}]$ is
  vertical, then $\rho$ is a \emph{regular singularity}.  Points where
  $\dim{\V_{\rho}[\Rc{q}]}>1$ are called \emph{irregular singularities}.
\end{definition}

It is not difficult to show that the regular points form a dense subset and
at them the classical existence and uniqueness theorems for differential
equations hold (see e.\,g.\ the discussion in \cite{wms:aims}).  At
singularities, one typically looses in particular the uniqueness -- there
can be anything from two to infinitely many one-sided solutions -- and the
existence of two-sided solutions going \emph{through} the point can become
a highly non-trivial problem.  Conventional numerical integrators will
break down when approaching such points.  For an analysis of the solution
behaviour in their neighbourhood, we need the following, more general
notions of solutions.

\begin{definition}\label{def:gensol}
  A \emph{generalised solution} of the differential equation
  $\mathcal{R}_{q}$ is an integral curve
  $\mathcal{N}\subseteq\mathcal{R}_{q}$ of the Vessiot distribution
  $\mathcal{V}[\mathcal{R}_{q}]$, i.\,e.\ a one-dimensional submanifold
  such that
  $T_{\rho}\mathcal{N}\subseteq\mathcal{V}_{\rho}[\mathcal{R}_{q}]$ at
  every point $\rho\in\mathcal{N}$.  The projection
  $\pi_{0}^{q}(\mathcal{N})\subset\E$ is called a \emph{geometric
    solution}.
\end{definition}

Generalised solutions live in the jet bundle $\je{q}$ and not in the base
manifold $\E$.  If a section $\sigma$ defines a strong solution in the
sense of Definition \ref{def:diffeq}, then
$\im{j_{q}\sigma}\subseteq\Rc{q}$ is a generalised solution and
$\im{\sigma}\subset\E$ the corresponding geometric solution.  However, not
all geometric solutions are graphs of functions.  In fact, they are not
even necessarily smooth curves, as they arise via a projection.

Away from the irregular singularities, the Vessiot distribution can be
locally generated by a vector field.  In \cite{wms:aims}, it is shown that
the irregular singularities are stationary points of this vector field.
Thus the analysis of an implicit differential equation satisfying our
assumptions can be reduced to the study of an autonomous dynamical system
on the submanifold $\Rc{q}\subset\je{q}$.  In the next section we will
discuss how this idea can be realised numerically.

In this approach, regular singularities play no role at all.  There goes a
unique generalised solution through any regular singularity and the
singularity is a smooth point of this curve \cite{wms:aims}.  The singular
character of such a point becomes apparent only, if one tries to interpret
the obtained curve as the prolongation of the graph of a function, as this
is usually not possible.  Indeed, if one considers the associated geometric
solution, it usually exhibits a cusp underneath the regular singularity.
However, this does not affect in the least the numerical determination of
the generalised solution.

\section{Numerical Integration of Implicit Differential Equations}
\label{sec:numint}

Classical approaches essentially transform numerically an implicit equation
into either an explicit differential equation or a differential algebraic
equation which is then solved by standard methods.  Such approaches run
into difficulties whenever the integration gets close to a singularity, as
at singularities the ranks of certain crucial matrices jump and thus
already in their vicinity condition numbers deterioriate.  Furthermore, at
singularities solutions are no longer unique and the number of solutions
can be anything from two to infinity.

We describe now a realisation of a different approach, namely the numerical
integration of the Vessiot distribution or more precisely, of a vector
field locally generating it (essentially the same approach was already used
by Tuomela \cite{jt:sing,jt:resol}).  Thus we determine directly
generalised solutions.  This approach has some advantages.  In particular,
it immediately resolves the problem of regular singularities, as we already
discussed above.  Equations of arbitrary high order can be tackled directly
without the need to rewrite them as first-order equations (in fact, in this
approach such a rewriting would lead to unnecessary large equations).
Finally, we will show that also the analysis of irregular singularities can
be supported by standard methods from dynamical systems theory, although
this remains a hard problem for larger systems.

The price to pay for these advantages is an increase in the system size.
If the original implicit system involves $m$ unknown functions and $k$
equations of maximal order $q$, then we have to deal with an autonomous
vector field in $(q+1)m+1$ unknown functions plus $k$ weak invariants of
this field.  For moderate values of $m$ and $q$, which we will exclusively
consider in this article, this increase is easily tolerable in view of the
advantages.  For larger systems, one could think of a combination of the
classical approach (away from the singularities) and the here described
approach (in the vicinity of singularities) to improve efficiency.

Let the differential equation $\Rc{q}\subset\je{q}$ be described by the
implicit system $\Fv(t,\uv_{(q)})=0$.  We want to approximate the unique
generalised solution starting at a given point $\rho^{(0)}\in\Rc{q}$ which
is not an irregular singularity.  We construct numerically a vector field
$X$ for which this solution is a trajectory and integrate it numerically
obtaining a sequence of points $\rho^{(i)}$.  From the discussion in the
previous section, it is obvious that the value $X_{\rho}$ of $X$ at some
point $\rho$ on the generalised solution generates the Vessiot space
$\V_{\rho}[\Rc{q}]$.  Therefore, having arrived at some point $\rho^{(n)}$,
the first task for constructing $\rho^{(n+1)}$ consists of determining
$\mathcal{V}_{\rho^{(n)}}[\mathcal{R}_{q}]$.

Computing the Vessiot space $\mathcal{V}_{\rho^{(n)}}[\mathcal{R}_{q}]$ at
a point $\rho^{(n)}\in\mathcal{R}_{q}$ is straightforward and requires in
principle only linear algebra.  Any vector
$X_{\rho^{(n)}}\in\mathcal{V}_{\rho^{(n)}}[\mathcal{R}_{q}]$ lies in the
contact distribution $\mathcal{C}^{(q)}|_{\rho^{(n)}}$ and thus can be
written as a linear combination of the basic contact fields given in
(\ref{eq:ciq}):
$X_{\rho^{(n)}}=aC^{(q)}_{\mathrm{trans}}|_{\rho^{(n)}}+\bv\Cv^{(q)}|_{\rho^{(n)}}$.
On the other hand, $X_{\rho^{(n)}}$ must be tangent to $\mathcal{R}_{q}$.
Hence, $X_{\rho^{(n)}}$ must satisfy the equations
$d\Fv|_{\rho^{(n)}}(X_{\rho^{(n)}})=0$.  Evaluation of this condition
yields the following homogeneous linear system of equations for the
coefficients $a$, $\bv$:
\begin{equation}\label{eq:vessdist}
  C^{(q)}_{\mathrm{trans}}(\Fv)(\rho^{(n)})a+
  \Cv^{(q)}(\Fv)(\rho^{(n)})\bv=0\,.
\end{equation}
By studying the ranks of the full coefficient matrix of this system and of
the submatrix $\Cv^{(q)}(\Fv)(\rho^{(n)})$ one can easily detect whether
$\rho^{(n)}$ is a singularity and if yes, what kind (see e.\,g.\
\cite{wms:aims}).

At a given point $\rho^{(n)}\in\Rc{q}$, it is straightforward to determine
a basis of the solution space of (\ref{eq:vessdist}) with standard
numerical methods.  Away from the irregular singularities, this basis
consists of a single vector $\tilde{X}_{\rho^{(n)}}$ generating
$\V_{\rho^{(n)}}[\Rc{q}]$.  To enhance the stability of the numerical
integration of $X$, we normalise this vector to unit length and obtain the
desired vector $X_{\rho^{(n)}}$.  In our implementation, the norm condition
is actually added to the linear system \eqref{eq:vessdist}.  On one side,
this makes the system nonlinear, but on the other side it also makes it
square and thus easily solvable by a Newton method.  This idea also helps
to overcome another problem: if $X_{\rho^{(n)}}$ is a unit length solution
of \eqref{eq:vessdist}, then the same holds for $-X_{\rho^{(n)}}$.  Thus if
one is careless, it may happen that one integrates the same piece of the
trajectory back and forth.  For the solution of the square nonlinear
system, we always take $X_{\rho^{(n-1)}}$ as starting value for the Newton
iteration.  For reasonably small step sizes, $X_{\rho^{(n)}}$ and
$X_{\rho^{(n-1)}}$ will not differ much so that we will obtain rapid
convergence to the right vector.

\begin{remark}
  For smaller systems, the linear system (\ref{eq:vessdist}) could be
  tackled symbolically treating the point $\rho^{(n)}\in\Rc{q}$ as a
  parameter.  As the behaviour of the system generally depends on the point
  $\rho^{(n)}$, one faces the nontrivial problem of solving a parametric
  linear system with potentially many necessary case distinctions.  A
  method for this was presented by Sit \cite{sit:sla} and it is also
  possible to use a Thomas decomposition (see e.\,g.\ \cite{bglr:thomas}
  and references therein).  Both approaches require that the parameter
  dependency is polynomial and are computationally quite demanding.  We
  will use the second approach elsewhere for the development of an
  effective theory of algebraic differential equations.  In this work, we
  will restrict to a purely numerical approach.
\end{remark}

While the linear system (\ref{eq:vessdist}) can be written down for any
point $\rho\in\je{q}$, it is actually defined only on the submanifold
$\Rc{q}\subset\je{q}$.  If a parametrisation of this submanifold is known,
then one could rewrite the system and the vector field it describes in
these parameters.  However, in practise such a parametrisation is often not
available.  Therefore we prefer to work with the jet coordinates on
$\je{q}$ and thus with redundant coordinates.  During the numerical
integration, we must then ensure that our approximate solution stays on the
submanifold $\Rc{q}\subset\je{q}$.

More precisely, our approach leads to the numerical integration of the
autonomous system
\begin{equation}\label{eq:dynsys}
  t'=a(t,\uv_{(q)}),\ \
  \uv'=a(t,\uv_{(q)})\dot{\uv},\ \ \dots\ \
  \uv'_{i}=a(t,\uv_{(q)})\uv_{i+1},\ \ \dots\ \
  \uv_{q}'=\bv(t,\uv_{(q)})
\end{equation}
where $'$ denotes the derivative with respect to some variable $x$
parametrising our generalised solution, all variables $t,\uv_{(q)}$ are
considered as independent algebraic variables and the functions $a$, $\bv$
arise from solving the linear system \eqref{eq:vessdist}.  However, we are
only interested in solutions of \eqref{eq:dynsys} which lie on the manifold
$\Rc{q}$, i.\,e.\ we have the additional algebraic constraints
$\Fv(t,\uv_{(q)})=0$.  By construction, these constraints represent weak
invariants of \eqref{eq:dynsys}, since the vector field $X$ is everywhere
tangential to the manifold $\Rc{q}$.

Hence we combined a standard numerical integrator applied to
\eqref{eq:dynsys} with a subsequent projection on the manifold $\Rc{q}$ (if
the obtained point $\rho^{(n+1)}$ lies too far away).  Such projections are
well-studied in numerical analysis and it is well-known that they do not
affect the convergence order of the numerical integrator.  Note furthermore
that solving \eqref{eq:vessdist} at a point $\rho$ satisfying
$\Fv(\rho)=\epsilon$ will lead to a vector which is tangential to the
manifold described by the perturbed system $\Fv(t,\uv^{(q)})=\epsilon$.
This observation implies that numerical errors will not lead to a strong
drift off $\Rc{q}$.  In fact, in many situations the manifold $\Rc{q}$ will
even be orbitally stable for the dynamical system \eqref{eq:dynsys}.
Indeed, we never observed any numerical problems due to a drift.

\section{Solutions at Irregular Singularities}
\label{sec:irreg}

As already mentioned above, irregular singularities become stationary
points of the vector field $X$ describing the Vessiot distribution locally
in their neighbourhood.  Generalised solutions through the singularity may
then be interpreted as one-dimensional invariant manifolds of this vector
field containing the singularity.  Hence for the analysis of the local
solution behaviour near an irregular singularity it is useful to compute
the invariant manifolds at the stationary point.

A more or less automated complete analysis of an irregular singularity is
in general only possible in certain low-dimensional situations.  If the
singularity corresponds to an hyperbolic stationary point, then the
dimension plays no role.  A complete analysis of this case will appear in
\cite{wms:quasilin}.  However, singularities correspond rarely to
hyperbolic stationary points (in the case of scalar higher-order equations
it is even impossible that a hyperbolic stationary point arises).  If the
centre manifold at a non-hyperbolic stationary point is at most
two-dimensional, then a complete analysis is possible using first a centre
manifold reduction to a two-dimensional system and then blow-ups
\cite{dla:qualplan} (and there even exist computer programmes for this task
like P4 described in \cite{dla:qualplan}).

The simplest situation arises, if an (un)stable or centre manifold is
one-dimensional.  In this case, it can immediately be identified with a
generalised solution through the stationary point.  If an invariant
manifold is higher dimensional, then one must analyse in more details the
reduced dynamics on it.  The key question is whether or not it is possible
to combine two trajectories approaching the singularity with the same
tangent to a smooth invariant manifold.  Assume for example that we are
dealing with a two-dimensional invariant manifold on which the phase
portrait looks like a node with two tangents.  Then almost all trajectories
will approach the singularity tangent to the eigenvector corresponding to
the eigenvalue whose real part has the smaller absolute value.  If this
eigenvector is transversal to the fibration $\pi^{q}$, then we obtain
generalised solutions through the singularity which correspond to smooth
classical solutions.  Otherwise, the classical solutions are of finite
regularity.

\begin{remark}
  Centre manifolds lead to further challenges.  It is well known that a
  centre manifold is not necessarily unique.  In fact, if we one has
  e.\,g.\ a saddle node, then there is a unique centre manifold on one side
  of the singularity and infinitely many centre manifolds on the other
  side.  In dynamical systems theory, one usually considers different
  centre manifolds as equivalent, as they are exponentially close when
  approach the singularity.  For us, each centre manifold corresponds to a
  different generalised solution and all of them are of interest.  Thus the
  analysis of the singularity requires in such a case precise statements
  about the (non-)uniqueness and the regularity of the centre manifolds
  which are not so easy to obtain.
\end{remark}

A complete analysis is always possible for a first-order scalar equation
$\Rc{1}\subset\je{1}$.  In this case, the Vessiot distribution induces a
dynamical system on the two-dimensional manifold $\Rc{1}$ and all its
stationary points can be studied using the methods described in
\cite{dla:qualplan}.  A concrete example is considered in Section
\ref{sec:examples}.  For higher-order scalar equations
$\Rc{q}\subset\je{q}$, one obtains a dynamical system on a
$(q+1)$-dimensional manifold and it is easy to see that no stationary point
of it can be hyperbolic, as the $q-1$ ``middle rows'' of the Jacobian are
multiples of the first row.  Hence we obtain many zero eigenvalues.

\begin{remark}\label{rem:quasilin}
  It was shown in \cite{wms:singbif} that quasi-linear
  equations\footnote{It is well-known that the fibration
    $\pi^{q}_{q-1}:\je{q}\rightarrow\je{q-1}$ defines an affine bundle.  A
    differential equation $\Rc{q}\subset\je{q}$ is \emph{quasi-linear}, if
    it is an affine subbundle.} have their own theory (see also the
  forthcoming work \cite{wms:quasilin} for a much more extensive treatment
  of this special case).  Given a quasi-linear equation
  $\Rc{q}\subset\je{q}$, one may consider instead of the Vessiot
  distribution on $\Rc{q}$ its well-defined projection into $\je{q-1}$.
  Furthermore, this projection is usually extendable to the whole jet
  bundle $\je{q-1}$ which leads to further phenomena specific to
  quasi-linear equations.  This projectability simplifies the analysis in
  the sense that there is no longer the need to work with redundant
  coordinates in an ambient space.  Otherwise, the analysis is performed
  along the same lines as for fully nonlinear systems by studying the
  stationary points of the projected Vessiot distribution which we call
  \emph{impasse points} to distinguish them from the singularities living
  one order higher.  A concrete example of a quasi-linear second-order
  equation will be studied in Section \ref{sec:examples}.
\end{remark}

For the computation of invariant manifolds, we implemented an algorithm for
the construction of a Taylor series approximation of the invariant manifold
originally developed by Beyn and Kle\ss~\cite{bk:invman} and later improved
by Eirola and von Pfaler \cite{ep:teim}.  Consider an $n$-dimensional
autonomous dynamical system $\dot{\xv}=\fv(\xv)$ with a stationary point
$\xi\in\RR^{n}$.  Assume that the spectrum of the Jacobian $J(\xi)$ of
$\fv$ in $\xi$ can be disjointly split into two parts,
$\mathrm{Spec}\bigl(J(\xi)\bigr)=\Sigma\sqcup\tilde{\Sigma}$, with a
corresponding splitting into generalised eigenspaces
$\RR^{n}=E\oplus\tilde{E}$.  It is well known that then, under certain gap
conditions, an invariant manifold $W$ exists which is tangent to $E$ and
which can locally be written as a graph over some neighbourhood $U$ of the
origin in $E$: $W=\bigl\{(\yv,\hv(\yv))\mid \yv\in U\subseteq E\bigr\}$.
Furthermore, the reduced dynamics on $W$ is given by an autonomous system
$\dot{\yv}=\gv(\yv)$ of dimension $\dim{E}$.  The goal is the construction
of a Taylor series approximation of $\hv$ and consequently of $\gv$.  Beyn
and Kle\ss~\cite{bk:invman} showed how this problem can be reduced to
solving multilinear Sylvester equations; in the subsequent improvement by
Eirola and von Pfaler \cite{ep:teim}, solving standard Sylvester equations
is sufficient.  Our implementation uses for this the built-in
\textsc{Matlab} procedure.

It should be noted that our implementation does not check whether
appropriate gap conditions are indeed satisfied and thus an invariant
manifold of sufficiently high regularity really exists.  This is the
responsibility of the user.  In our current implementation, it is only
possible to determine Taylor series up to degree 10.  The reason is simply
that certain combinatorial coefficients always appearing in the
computations independent of the concrete dynamical system considered have
been precomputed and stored on file.  So far, this precomputation has been
done only up to degree 10, but an extension to higher degree would be
possible without problems.

\begin{remark}
  One may consider this part of our work as a combined numerical-symbolic
  computation.  While the actual computation of the invariant manifold $W$
  is done purely numerically, its output -- Taylor polynomials for $\gv$
  and $\hv$ -- is in symbolic form.  Indeed, in certain situations, e.\,g.\
  for visualisations of the reduced dynamics, we will use the output as
  symbolic input for further computations.  A concrete example will appear
  in Section~\ref{sec:examples}.
\end{remark}

\section{Visualisation of Implicit Differential Equations}
\label{sec:visual}

In low-dimensional situations, a visualisation of the generalised solutions
is very useful for an understanding of the solution behaviour of an
implicit differential equation.  One fundamental problem of such a
visualisation is to obtain evenly spaced streamlines filling the whole area
of interest.  A number of solutions have been developed for planar vector
fields.  We have chosen to follow the approach of Jobard and Lefer
\cite{jl:stream} (which also underlies the \textsc{MuPAD}
\texttt{streamlines} command) and to adapt it for our purposes.

There are three different situations where a more or less complete
visualisation is possible.  The first case concerns a scalar first-order
equation $\Rc{1}\subset\je{1}$.  Although we are then dealing with a
two-dimensional vector field, it does not live on the plane $\RR^{2}$ (as
always assumed by Jobard and Lefer), but on the two-dimensional submanifold
$\Rc{1}$ lying in a three-dimensional ambient space $\je{1}\cong\RR^{3}$
(for the trivial fibration $\pi:\RR\times\RR\rightarrow\RR$).  In the
literature, this is usually called a \emph{2.5D visualisation}.

As already mentioned in Remark \ref{rem:quasilin}, for quasi-linear
equations, a reduction is possible which allows for a visualisation of
slightly higher-dimensional situations.  For a scalar quasi-linear
first-order equation, we obtain this way a two-dimensional planar vector
field and thus can use a standard \emph{2D visualisation}.  For a scalar
quasi-linear second-order equation or for a system of two first-order
quasi-linear equations, the projection yields a three-dimensional vector
field on $\RR^{3}$ and thus leads to a classical \emph{3D visualisation}.

Our implementation covers all three cases: 2D, 2.5D and 3D visualisation.
In the 2D case, we can essentially use the basic form of the algorithm of
Jobard and Lefer with only one minor modification.  We briefly recall its
basic ideas and refer for more details to their original paper
\cite{jl:stream}.  The algorithm works with \emph{seed points} which are
used as initial data for the computation of trajectories (in both
directions).  Besides the vector field, it is given as input a first seed
point which can essentially be chosen randomly and two parameter
$0<d_{\mathrm{test}}<d_{\mathrm{sep}}$ prescribing the desired minimal and
average distance, resp., between two streamlines.

As first step, the two semitrajectories starting at the initial seed point
are computed until they reach the boundary of the plotting region which
gives the first streamline.  As this computation is done numerically, the
streamline actually consists of a list of \emph{sample points}.  We produce
new potential seed points by going from each sample point orthogonally to
the streamline the distance $d_{\mathrm{sep}}$ (in positive and negative
direction).  For each potential seed point, it must then be checked that it
is still in the plotting region and that its distance from any sample point
is at least $d_{\mathrm{sep}}$.  If not, the point is discarded.  Now a new
seed point is picked randomly and the corresponding streamline computed
until it either reaches the boundary of the plotting region or it gets
closer than $d_{\mathrm{test}}$ to some already computed sample point.
When all sample points on the new streamline are obtained, all currently
collected seed points must be checked again whether they have a distance
greater than $d_{\mathrm{sep}}$ from them.  The whole process is iterated,
until no admissible seed points exist any more.

The one above mentioned modification concerns the treatment of impasse
points.  For a quasi-linear first-order equation, the impasse points are
the stationary points of the vector fields whose streamlines we want to
compute.  Jobard and Lefer do not mention any special treatment of such
points.  However, it is well known that in particular in the neighbourhood
of non-hyperbolic stationary points the dynamics can be quite complicated
and difficult to resolve numerically (e.\,g.\ if elliptic sectors exist).
To avoid possible numerical problems, our implementation takes as further
input a list of the impasse points (if there is a whole curve of impasse
points, then it is represented by a list of sufficiently close sample
points) which are considered as additional sample points corresponding to
degenerate streamlines and a parameter $d_{\mathrm{s}}$ prescribing the
minimal distance from an impasse point.  In practise, one tries to choose
$d_{\mathrm{s}}$ as small as possible without encountering numerical
problems, as the neighbourhood of an impasse point is of course of
particular interest.

Computationally, the most expensive part of this algorithm is not the
numerical integration but the many checks whether points are sufficiently
far away from the already computed sample points.  As the number of sample
points increases with every additional streamline, this process becomes
more and more expensive.  As an optimisation, the plotting region is
divided into squares and a point in one square is only compared with sample
points in the same square or neighbouring squares.  Furthermore, the
numerical integration must be adapted to these tests.  In our
implementation, we work for simplicity with a constant step size $h$ and
take as a new sample point the point obtained after
$\lceil h/d_{\mathrm{sep}}\rceil$ integration steps.  More refined versions
using e.\,g.\ continuous Runge-Kutta methods with variable step sizes are
possible and probably necessary for stiff vector fields, but in our
experiments our simpler approach always produced good results.

In the 2.5D case, we are still on a 2D manifold, but it lives in an ambient
three-dimensional space.  This requires some adaptions of the above
described approach.  For producing new seed points, we must choose a
direction orthogonal to the streamline \emph{and} tangential to the
manifold in order to obtain a unique direction.  Nevertheless, the thus
obtained point will generally not lie on the manifold and must be
orthogonally projected back to the manifold.  However, the projection
changes the distance from the starting point.  Hence we must walk for a yet
undetermined distance and then project which leads to a non-linear system
of equations for two parameters.  Furthermore, one must discuss what
``distance'' actually should mean.  We use the 3D Euclidean distance
instead of some intrinsic distance within the manifold which is
computationally much simpler and turned out to be sufficient for all
practical purposes.  As a side effect, we divide now the (3D) plotting
range into cubes for optimising the distance tests.  As in the 2D case, the
user must provide a list with all irregular singularities of the considered
differential equation and a distance parameter $d_{\mathrm{s}}$.

In the 3D case, we must mainly adapt the production of potential seed
points.  Opposed to the two cases treated so far, there is no distinguished
direction to walk from an already computed sample point.  Instead, we
consider now a circle around it with radius $d_{\mathrm{sep}}$ in the plane
orthogonal to the current streamline.  Elementary geometric computations
show that we can place six points on the circle such that no two of them
are closer than $d_{\mathrm{sep}}$.  Thus each sample point yields six new
potential seed points.  Otherwise, we proceed exactly as in the 2D case.

The thus adapted algorithm produces now streamlines which are evenly spaced
in 3D.  Unfortunately, this is generally not sufficient for producing good
pictures.  The problem is that a 3D picture is projected onto some image
plane and the projected streamlines will no longer be evenly spaced.  In
fact, projected streamlines will cover each other or intersect which makes
it quite hard to interpret the obtained images.  This effect is often
called \emph{visual clutter}.  Our programme tries to enhance the 3D
visibility by a postprocessing step in which colour, transparency and
thickness of the streamlines are modified according to their position
relative to the observer.  However, as one can see in a concrete example in
Section \ref{sec:examples} below, the effect is limited.

Further improvements can probably be achieved by implementing further
visualisation algorithms specifically designed for 3D vector fields.  The
literature provides a number of such algorithms (see e.\,g.\
\cite{cyymm:3dvf,mchm:3dvf}), but it is unclear which one is best suited
for our application.  Furthermore, we will see in Section
\ref{sec:examples} below that in many situations one has here not only a
visualisation problem, but actually also mathematical problems.  Another
alternative would be the use of more specialised 3D rendering software like
\textsc{ParaView}\footnote{\url{https://www.paraview.org}} allowing for
many special effects.  In particular, some of these programmes allow to
``fly into the 3D image''.  Our evenly spaced streamlines should represent
an ideal starting point for such a presentation.

\section{Examples}
\label{sec:examples}

We discuss now the use of our \textsc{Matlab} programmes in the analysis of
two concrete implicit differential equations.  Both examples stem from
actual applications and their singularities will be discussed in more
details elsewhere.  Here, we mainly present some visualisations produced
with our programmes and discuss some problems and shortcomings.

\subsection{A Scalar First-Order Equation}

In the context of reconstructing the position and orientation of known 3D
objects from a 2D image of them, the following fully non-linear scalar
first-order equation arises:
\begin{equation}\label{eq:3D2D}
  (1+t^{2})\dot{u}^{2} + u^{2} = r(t)^{2}
\end{equation}
where the function $r:\RR\rightarrow\RR_{>0}$ encodes certain information
about the 3D object.  Following our geometric approach, we consider
\eqref{eq:3D2D} as the description of a two-dimensional surface $\Rc{1}$ in
the three-dimensional jet bundle $\je{1}$.  The linear system
\eqref{eq:vessdist} determining the Vessiot space $\V_{\rho}[\Rc{1}]$ at a
point $\rho=(t,u,\dot{u})\in\Rc{1}$ reduces here to the single equation
\begin{equation}\label{eq:VD3D2D}
  (t\dot{u}^{2}-r(t)\dot{r}(t)+u\dot{u})a+(1+t^{2})\dot{u}b=0\,.
\end{equation}
We see that any point with $\dot{u}=0$ is a singularity, as at such points
the coefficient of $b$ vanishes implying that either $a=0$ (regular
singularity) or $\dim{\V_{\rho}[\Rc{1}]}=2$ (irregular singularity).  Thus
the singularities define two curves $(t,\pm r(t),0)$ with $t\in\RR$.  Since
we assume that always $r(t)>0$, irregular singularities are characterised
by the additional condition $\dot{r}(t)=0$ (which implies that
\eqref{eq:VD3D2D} reduces to $0=0$), i.\,e.\ they correspond to the
critical points of these curves.  Eq.~\eqref{eq:VD3D2D} is easily solved
away from the irregular singularities and one finds that the Vessiot
distribution is almost everywhere generated by the vector field
\begin{displaymath}
  X=(1+t^{2})\dot{u}(\partial_{x}+\dot{u}\partial_{u})+
       (r(t)\dot{r}(t)-u\dot{u}-t\dot{u}^{2})\partial_{\dot{u}}\,.
\end{displaymath}

\begin{figure}[ht]
  \centering
  \includegraphics[height=6.7cm]{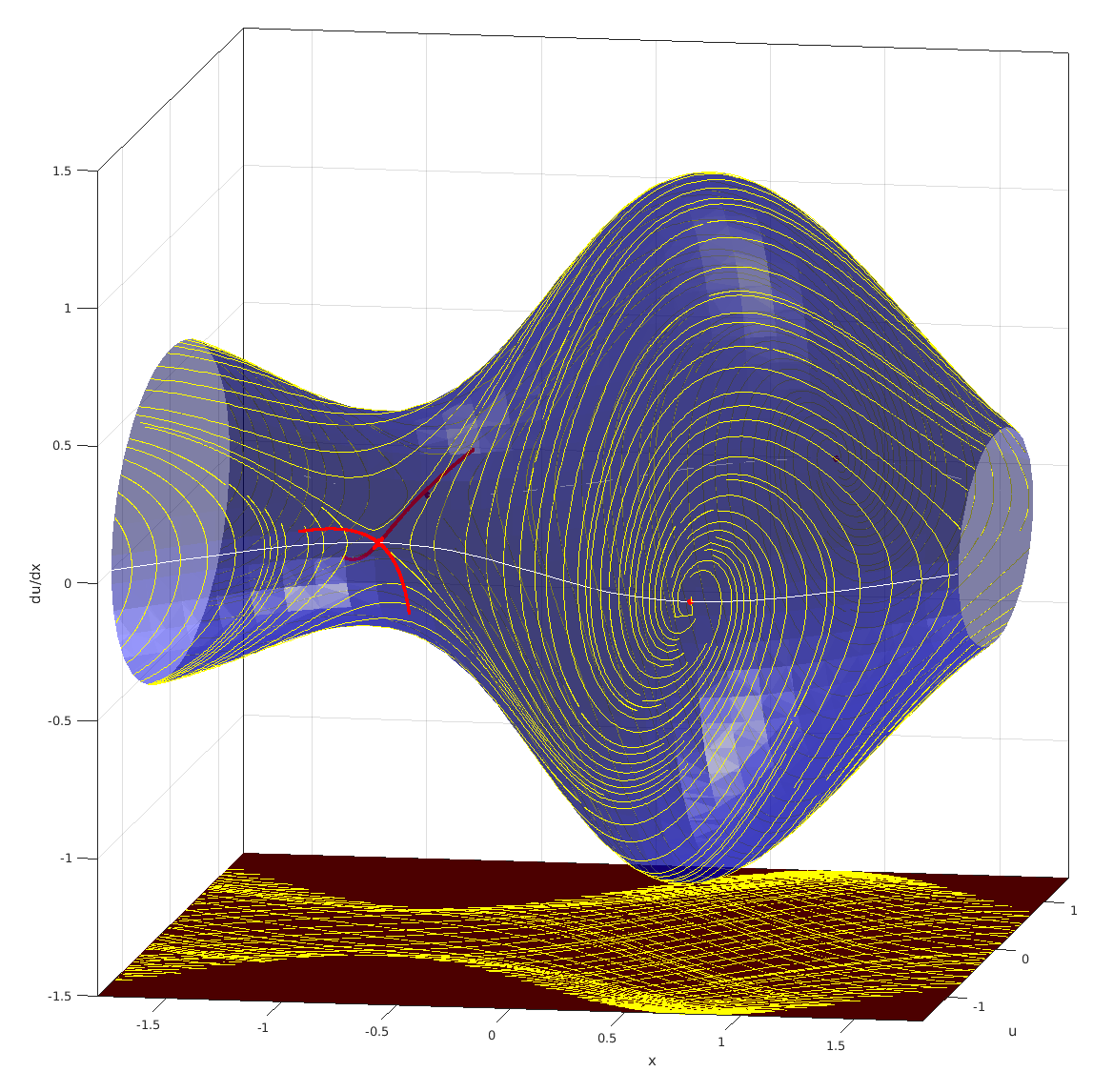}\quad
  \includegraphics[height=6.7cm]{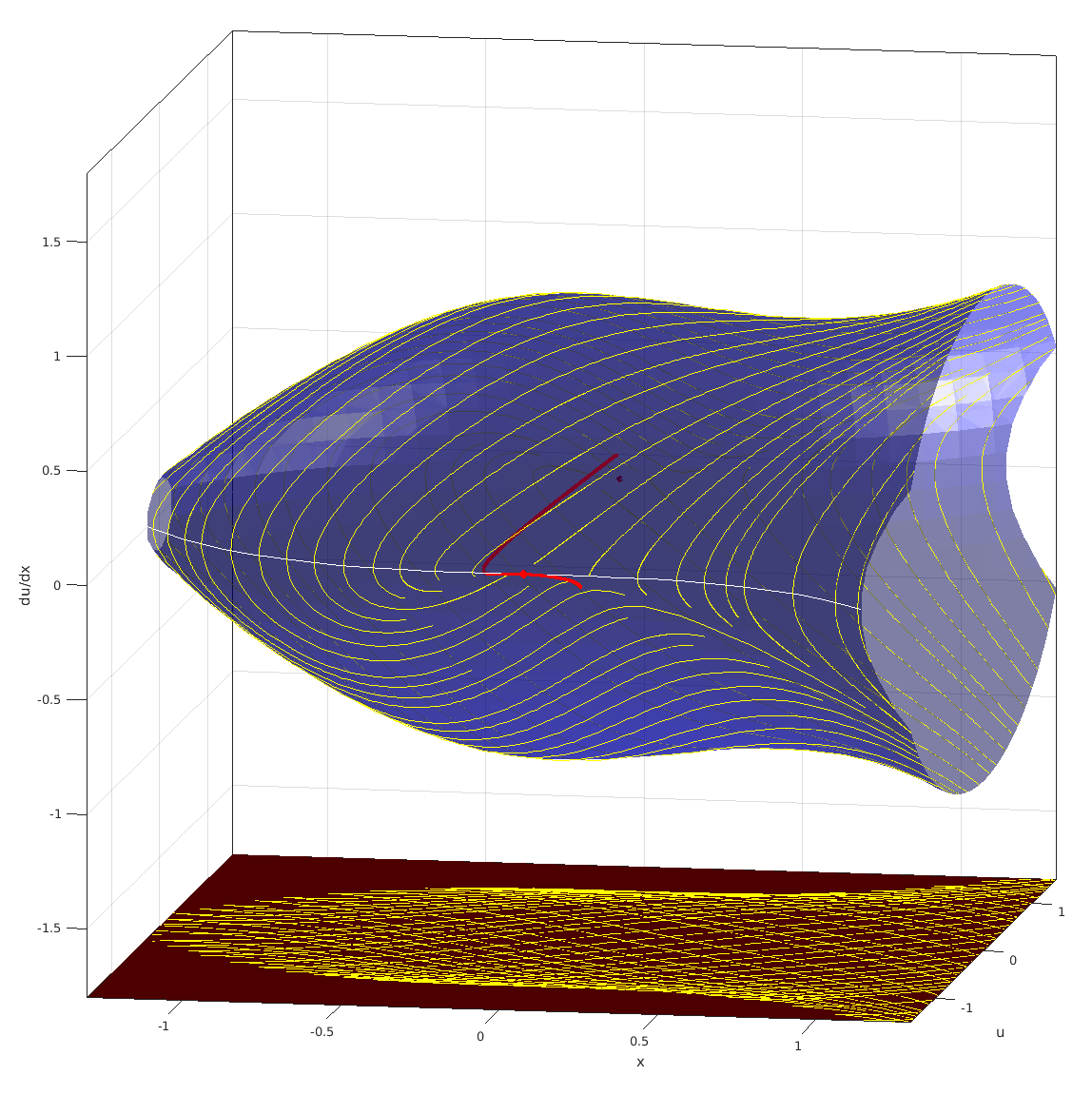}
  \caption{$\Rc{1}\subset\je{1}$ given by \eqref{eq:3D2D} with generalised
    and geometric solutions. On the left hand side for the choice
    $r(t)=1+\sin{(2t)}/2$; on the right hand side for $r(t)=1+t^{3}/3$}
  \label{fig:3D2D-1}
\end{figure}

Figure \ref{fig:3D2D-1} shows the surface $\Rc{1}$ for two choices of the
function $r(t)$.  The white line depicts the singularities; irregular ones
are marked as red points.  The yellow lines on the surface represent
generalised solutions computed with the 2.5D version of the algorithm by
Jobard and Lefer.  At the bottom, one can see the corresponding geometric
solutions obtained by a simple projection.  Because of their large number,
it is not easily visible in the projection that they change direction
whenever the generalised solution crosses the white line; but one clearly
recognises this behaviour from the form of the generalised solutions.

For a closer analysis of the behaviour at an irregular singularity
$\rho=(\bar{t},\pm r(\bar{t}),0)$, one needs the Jacobian of the vector
field $X$ at $\rho$.  It is easily determined:
\begin{displaymath}
  J=
  \begin{pmatrix}
    0 & 0 & (1+\bar{t}^{2})\\
    0 & 0 & 0 \\
    r(\bar{t})\ddot{r}(\bar{t}) & 0 & r(\bar{t})
  \end{pmatrix}\,.
\end{displaymath}
One of its eigenvalues is always $0$ with eigenvector
$\left(\begin{smallmatrix} 0\\ 1\\ 0 \end{smallmatrix}\right)$.  However,
as this eigenvector is not tangential to $\Rc{1}$, we must discard this
eigenvalue.\footnote{Note that, strictly speaking, we are dealing here with
  a vector field on a two-dimensional manifold.  If we had a nice
  parametrisation of the manifold, we could express the vector field $X$ in
  these parameters and would obtain a $2\times2$ Jacobian.  As it is in
  general difficult to find such parametrisations, we use instead the three
  coordinates of the ambient space $\je{1}$.  Consequently, we obtain a too
  large Jacobian and must see which two eigenvalues are the right ones.
  This is easily decided by checking whether the corresponding
  (generalised) eigenvectors are tangential to $\Rc{1}$.} If we write
$\beta=r(\bar{t})^{2}+4r(\bar{t})\ddot{r}(\bar{t})(1+\bar{t}^{2})$, then
the other two -- relevant -- eigenvalues are given by
$(r(\bar{t})\pm\sqrt{\beta})/2$.  There arises now a total of five
different cases depending on the value of $\beta$.  Three of them can be
seen in Fig.~\ref{fig:3D2D-1}.  In the left picture, one can see a saddle
point (arising for $\beta>r(\bar{t})^{2}$) and a focus (arising for
$\beta<0$).  The right picture shows the case of a semi-hyperbolic
stationary point with a one-dimensional centre manifold (arising for
$\beta=r(\bar{t})^{2}$ which is equivalent to $\ddot{r}(\bar{t})=0$).  In
this case, further subcases have to be distinguished depending on the
values of higher derivatives of $r$ at $\bar{t}$.

Fig.~\ref{fig:3D2D-1} also shows (some of) the invariant manifolds at the
irregular singularities.  At the saddle point on the left hand side, one
can see the stable and the unstable manifold as red lines (the reader may
choose which one is the stable manifold, as one could perform the same
analysis with the vector field $-X$ for which the eigenvalues just swap
signs).  As each of these manifolds defines a generalised solution going
through the irregular singularity, we conclude that here two generalised
solutions intersect.  At the focus, no (real) invariant manifolds exist.
The generalised solutions approach the irregular singularity
asymptotically, however without a well-defined tangent.  Hence, here we
have no generalised solution \emph{through} the singularity.

In the right picture, a short piece of (an approximation of) a centre
manifold is shown.  A closer analysis of the behaviour at this irregular
singularity (which is beyond the scope of this article) reveals that we are
here actually in a situation where no analytic centre manifold exists (in
fact, where the centre manifolds are probably only of finite regularity),
i.\,e.\ where the Taylor series approximations computed by our algorithm do
not converge.  This behaviour leads to numerical problems for the algorithm
and in such cases one can often determine only experimentally a reasonable
order of approximation where the computation succeeds and produces a
reasonable result.  In our case one can see that the computed approximation
probably describes qualitatively correctly the form of the centre manifold
but that the piece shown can be accurate only rather close to the
singularity, as further away it intersects with other generalised solutions
which is not possible.  The shown red line approximates one generalised
solution going through the irregular singularity.  Generally, centre
manifolds are not unique and then each centre manifold yields a different
generalised solution. It is a classical result that all different centre
manifolds are exponentially close to each other and thus possess the same
Taylor polynomial (to any order for which it exists).  As our algorithm is
based on computing a Taylor series approximation, it cannot distinguish
different centre manifolds.

\subsection{A Quasi-Linear Second-Order Equation}

In an optimal control problem related to financial economics, the following
scalar quasi-linear second-order equation arises \cite{bcw:singde}
\begin{equation}\label{eq:bcw}
  t^{2}\ddot{u}=at\dot{u}+bu-c(\dot{u}-1)^{2}
\end{equation}
together with the initial conditions $u(0)=0$ and $\dot{u}(0)=1$.  Here
$a,b,c\in\RR$ are parameters.  To avoid case distinctions, we assume that
$bc\neq0$ and $a+b\neq0$.  Following the above mentioned reduction process
for quasi-linear equations, we are lead to study the three-dimensional
vector field
\begin{equation}\label{eq:Ybcw}
  Y=t^{2}\partial_{t} + t^{2}\dot{u}\partial_{u} +
       \bigl(at\dot{u}+bu-c(\dot{u}-1)^{2}\bigr)\partial_{\dot{u}}\,.
\end{equation}
Under the above made assumptions on the parameters, the stationary points
of this vector field -- and thus the impasse points of \eqref{eq:bcw} --
lie on the parabola $t=0$ and $bu=c(\dot{u}-1)^{2}$.  At the ``tip'' of the
parabola, i.\,e.\ at the point $\rho=(0,0,1)$, the Jacobian of $Y$ has $0$
as a triple eigenvalue (and a non-trivial Jordan normal form with two
blocks).  At all other points on the parabola, the Jacobian has only a
double eigenvalue $0$.  It is still an open problem to study in detail the
local solution behaviour around $\rho$; at the other impasse points a
centre manifold reduction yields a planar problem which can be completely
analysed.

\begin{figure}[ht]
  \centering
  \includegraphics[height=6cm]{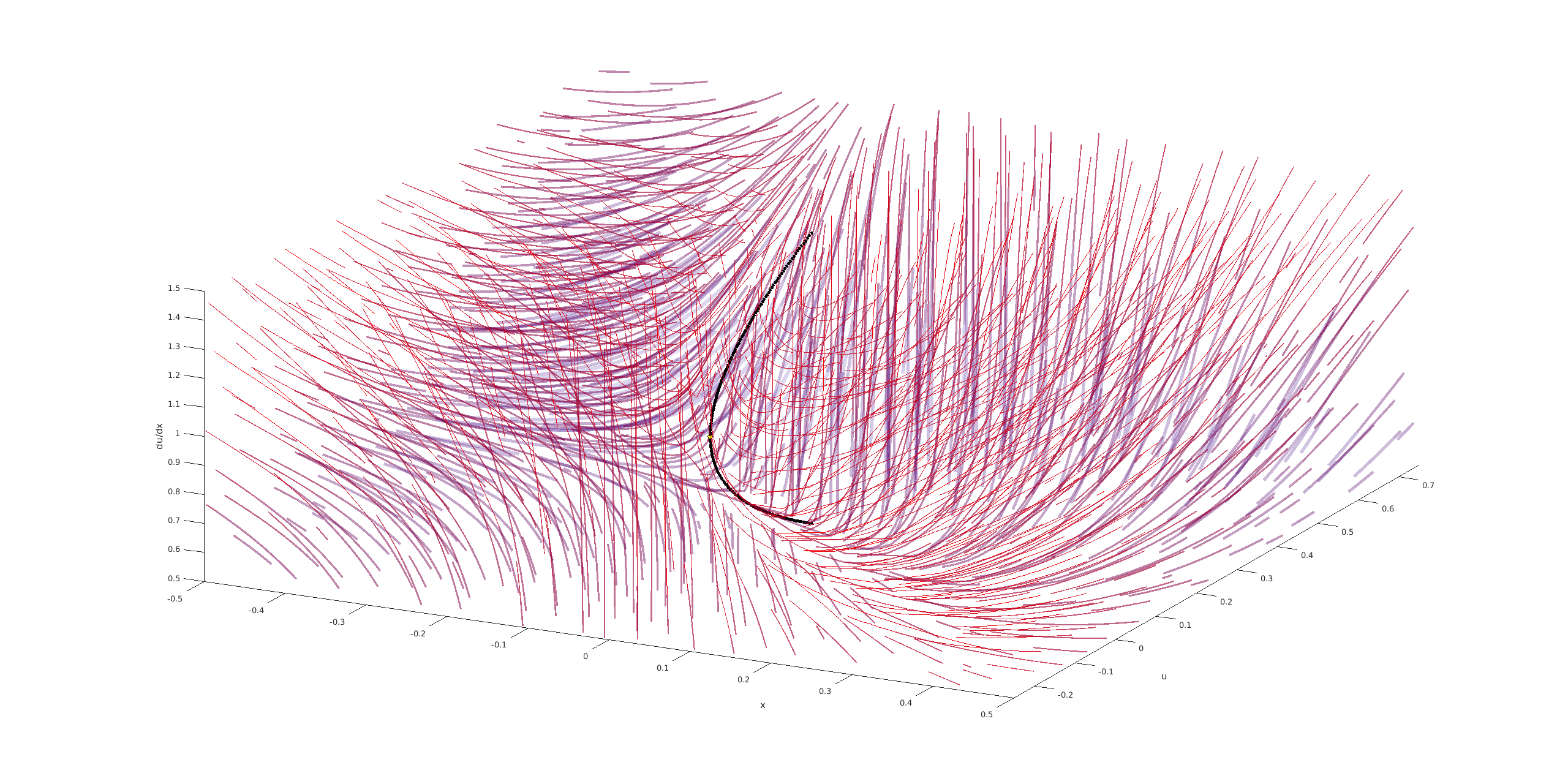}
  \caption{3D streamlines of the vector field $Y$ given by \eqref{eq:Ybcw}}
  \label{fig:wink3d}
\end{figure}

Fig.~\ref{fig:wink3d} shows a 3D visualisation of the vector field $Y$
given by \eqref{eq:Ybcw} using the parameter values $a=b=c=1$.\footnote{In
  \cite{bcw:singde} it is shown that for these parameter values there are
  infinitely many solutions reaching the ``tip'' $\rho$.}  The black
parabola contains the impasse points.  One clearly sees that the extension
of the algorithm by Jobard and Lefer to 3D only partly helps with the
visualisation.  The above mentioned problem that projected streamlines
cover each other or intersect is clearly visible and drastically reduces
the interpretability of the obtained images.  However, one should keep in
mind that our goal is to obtain a better understanding of the dynamics in a
neighbourhood of the ``tip'' $\rho$ of the parabola.  Here, the real
problem is less the 3D visualisation but a mathematical one.  If we get too
close to the singularity, the numerical integration will break down.  As
mentioned above, we introduced specifically for this purpose a parameter
$d_{\mathrm{s}}$ as stopping criterion for the numerical integration when
we approach a singularity.  In the example at hand, we could actually get
fairly close to the parabola and thus a careful study of the picture
reveals at least some indications about the local solution behaviour, in
particular if we combine it with the information obtained from the next
picture.

\begin{figure}[ht]
  \centering
  \includegraphics[height=6cm]{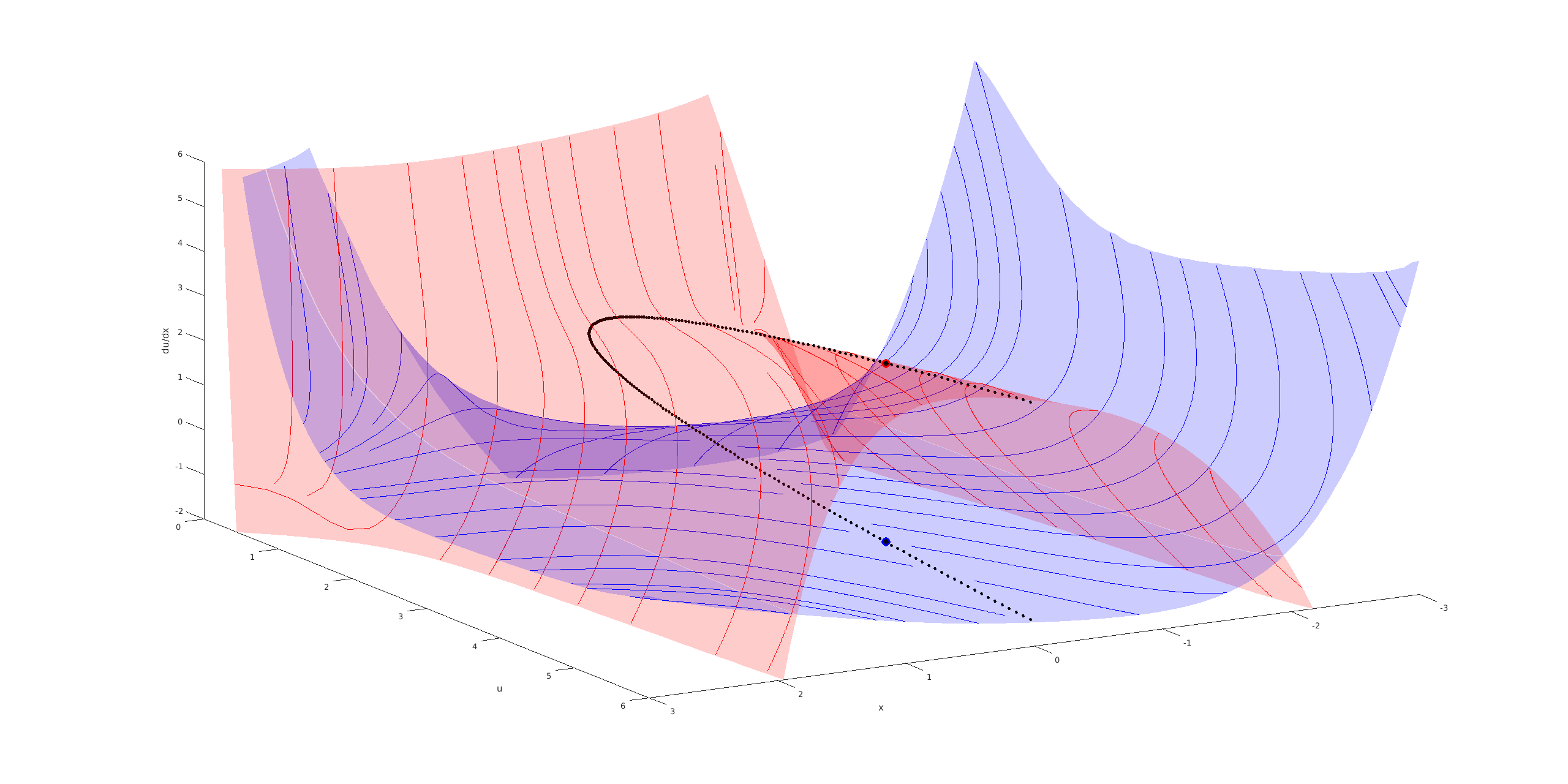}
  \caption{Centre manifolds and reduced dynamics on them for two different
    impasse points of \eqref{eq:bcw}}
  \label{fig:winkreddyn}
\end{figure}

Fig.~\ref{fig:winkreddyn} shows approximations of two-dimensional centre
manifolds at two points on the parabola away from the ``tip'' plotted
together with the reduced dynamics on them (one point and the corresponding
manifold are shown in red, the other one in blue).  As the parabola
consists entirely of stationary points, it is part of any centre manifold
through a point on it which is clearly visible in the picture.  Indeed, if
one computes the eigenvectors of the Jacobian of $Y$ to the eigenvalue $0$
(away from the ``tip'' $\rho$ the Jacobian is always diagonalisable), then
one of them is tangential to the parabola.  The other provides the
direction of the reduced dynamics which consists of streamlines cleanly
intersecting the parabola.

Fig.~\ref{fig:winkreddyn} was obtained by combining two of the above
described algorithms.  First, a Taylor series approximation of the centre
manifolds and the reduced dynamics on it are computed via the approach of
Eirola and von Pfaler.  Then, a 2.5D visualisation of the reduced dynamics
is determined following the method of Jobard and Lefer.  We have plotted
rather large parts of the centre manifold approximations and one probably
should analyse the quality of the approximation further away from the
selected impasse points in more details.  The picture seems to indicate
that the centre manifolds intersect.  While this is surely possible, the 3D
picture shown in Fig.~\ref{fig:wink3d} gives no indication that such a
phenomenon actually occurs.

\section{Conclusions}

In this article, we briefly sketched the use of geometric techniques, in
particular of the Vessiot distribution, for the analysis of implicit
ordinary differential equations.  We described a suite of \textsc{Matlab}
programmes supporting such an analysis with numerical integrations and
visualisations.  While there is surely still much room for improving the
programmes, in particular concerning the 3D visualisation, they already
proved useful for the analysis of concrete differential equations.

Our discussion of two concrete examples in Section \ref{sec:examples}
clearly shows that the analysis of singularities or impasse points cannot
be completely automatised even for low-dimensional problems.  One needs
further (mainly symbolic) computations, firstly to know what pictures could
be useful for the analysis and then for a better understanding of what the
picture are actually showing.  Nevertheless, in our experience pictures
like the ones shown in Section \ref{sec:examples} are very helpful for
unravelling the local solution behaviour around singularities.

\bibliography{NumVisImplODE.bib}
\bibliographystyle{plain}

\end{document}